# Stable determination of generic simple metrics from the hyperbolic Dirichlet-to-Neumann map

Plamen Stefanov*
Department of Mathematics
Purdue University
West Lafayette, IN 47907, USA

Gunther Uhlmann†
Department of Mathematics
University of Washington
Seattle, WA 98195, USA

**Abstract**

We prove Hölder type stability estimates near generic simple Riemannian metrics for the inverse problem of recovering such metrics from the Dirichlet-to-Neumann map associated to the wave equation for the Laplace-Beltrami operator.

## 1 Introduction and Main results

In this paper we consider the inverse problem of determining a Riemannian metric on a Riemannian manifold with boundary from the vibrations measured at the boundary. This information is encoded in the hyperbolic Dirichlet-to-Neumann (DN) map associated to the solutions to the wave equation. We concentrate on the stability question, that is if two hyperbolic DN maps are close in an appropriate topology, how close are the Riemannian metrics? We apply stability results obtained recently by the authors for the boundary rigidity problem [StU1], [StU2] to study this problem. We now describe in more detail the problem and the results.

Let $(M, g)$ be a compact Riemannian manifold with boundary. We denote by $\Delta_g$ the Laplace-Beltrami operator. In local coordinates $g(x) = (g_{ij}(x))$ it is given by

$$\Delta_g = (\det g)^{-\frac{1}{2}} \sum_{i,j=1}^{n} \frac{\partial}{\partial x^i} (\det g)^{\frac{1}{2}} g^{ij} \frac{\partial}{\partial x^j}.$$

in $\Omega$. Here $(g^{ij}) = (g_{ij})^{-1}$, $\det g = \det(g_{ij})$, and we will use freely the convention of raising and lowering indices of tensors. Consider the following problem

$$\begin{cases} (\partial_t^2 - \Delta_g)u = 0 & \text{in } (0, T) \times M, \\ u|_{t=0} = \partial_t u|_{t=0} = 0 & \text{in } \Omega, \\ u|_{(0,T) \times \partial \Omega} = f, \end{cases} \quad (1)$$

where $f \in C_0^1(\mathbf{R}_+ \times \partial \Omega)$. Denote by $\nu = \nu(x)$ the outer normal to $\partial M$ at $x \in \partial M$, normalized so that $g^{ij}\nu_i \nu_j = 1$. We define the hyperbolic Dirichlet-to-Neumann (DN) map $\Lambda_g$ by

$$\Lambda_g f := \sum_{i=1}^{n} \nu^i \frac{\partial u}{\partial x^i} \bigg|_{(0,T) \times \partial M},$$

*Partly supported by NSF Grant DMS-0400869
†Partly supported by NSF and a John Simon Guggenheim fellowship

where $v^i = \sum g^{ij}v_j$. It is easy to see that if

$$\psi : M \to M$$

is a diffeomorphism with $\psi|_{\partial M} = \text{Id}$, then $\Lambda_{\psi^*g} = \Lambda_g$, where $\psi^*g$ denotes the pull back of the metric g. The inverse problem is therefore formulated in the following way: knowing $\Lambda_g$, can one determine the metric $g$ up to an isometry that leaves the boundary fixed?

An affirmative answer to this question for smooth metrics was given by Belishev and Kurylev [BK]. Their approach is based on the boundary control method introduced by Belishev [B1] and uses in a very essential way a unique continuation principle proven by Tataru [T]. Because of the latter, it is unlikely that this method would prove Hölder type stability estimates even under geometric and topological restrictions. We also refer to [KKL1], [B2], [KKLM] and the references therein for more uniqueness results in this direction.

Hölder type of conditional stability estimate was proven by the authors in [StU3] for metrics close enough to the Euclidean one in $C^k$, $k \gg 1$ in three dimensions. Hölder type stability estimates were proven in [IS] and [Su] for the hyperbolic DN map associated to the Euclidean wave equation plus a potential.

The conditional type of the estimate, typical for such kind of inverse problems, is due to an additional a priori condition of boundedness of the $H^s$ norm of the metrics for some large $s$. It can be considered as a compactness condition in $H^s$ with smaller $s$. A well-known functional analysis argument shows that under such compactness condition, the map $\Lambda_g \to g$ (by identifying isometric metrics) must be continuous, once we know it is well-defined. This was exploited in [AKKLT] and minimal geometric conditions guaranteeing the compactness condition were established there in terms of bounds of certain geometric invariants, depending only on the second derivatives of the metric. The continuity of the map above however, does not give information about the type of a possible stability estimate, i.e. about the modulus of continuity of that map, see also [KKL2]. In this paper, we prove Hölder type of stability near generic *simple* metrics. A Riemannian manifold $(M, g)$ is simple, if $M$ is simply connected, $\partial M$ is strictly convex and any two point in $M$ can be connected by a single minimizing geodesic depending smoothly on them, see Definition 1 in next section.

Since simple manifolds are diffeomorphic to the unit ball in the Euclidean space, from now on, without loss of generality, we consider the case that $M = \bar{\Omega}$ where $\Omega$ is a bounded domain in the Euclidean space with smooth boundary.

As we mentioned above, we use recent result by the authors about the so-called boundary rigidity problem. The latter can be formulated as follows: Let $g$ be a simple metric in $\Omega$. Can we determine $g$, up to an isometry as above, from the knowledge of the distance function $\rho_g(x, y)$, known for all $x$, $y$ on the boundary $\partial \Omega$? The main result in [StU2] is that for $k \gg 1$, this is true for a dense open set $\mathcal{G}^k(\Omega)$ in $C^k(\bar{\Omega})$ of simple metrics, and $\mathcal{G}^k$ is defined as the set of those metrics, for which the linearized problem, integrals of 2-tensors along geodesics, is s-injective, see section 2. Moreover, $\mathcal{G}^k$ contains all real analytic simple metrics in $\Omega$.

The main result of this paper is the following.

**Theorem 1** *There exist $k > 0$, $0 < \mu < 1$, such that for any $g_0 \in \mathcal{G}^k$ and $T > \text{diam}_{g_0}(\Omega)$, $0 < \varepsilon < T - \text{diam}_{g_0}(\Omega)$, there exists $\varepsilon_0 > 0$, with the property that if*

$$\|g_m - g_0\|_{C(\bar{\Omega})} < \varepsilon_0, \quad \|g_m\|_{C^k(\bar{\Omega})} \leq A, \quad m = 1, 2,$$

*with some $A > 0$, then one can find a $C^3(\bar{\Omega})$ diffeomorphism $\psi : \bar{\Omega} \to \bar{\Omega}$, $\psi|_{\partial \Omega} = \text{Id}$, such that*

$$\|g_1 - \psi_*g_2\|_{C^2(\bar{\Omega})} \leq C \|\Lambda_{g_1} - \Lambda_{g_2}\|^\mu_{H^1_0([0,\varepsilon]\times\partial\Omega) \to L^2([0,T]\times\partial\Omega)}. \tag{2}$$

**Remark 1.** The condition about the closeness of the metrics to $g_0$ can, of course, be formulated in an invariant way: for some pull backs $\psi_i^* g_i$ of $g_i$, $i = 1, 2$, with $\psi_i$ as above, we require that $\psi_i^* g_i$ be $\varepsilon_0$-close to $g_0$. One can also study orbits of metrics under actions of such diffeomorphisms and express that condition as distance between the orbits of $g_1$, $g_2$. We remark also that the a-priori condition on boundedness of the $C^k$ norm of $g_i$, $i = 1, 2$, can be formulated invariantly in terms of a bound of the covariant derivatives of the curvature tensor as in [LSU].

**Remark 2.** One can generalize this result to lower order perturbation of $\Delta_g$. More precisely, consider $P = -\Delta_g + \sum b^j \partial/\partial x^j + q$, where $b = \{b^j\}$ and $q$ are complex-valued. In this case, $\Lambda$ is also preserved under the transformation $P \mapsto a^{-1} Pa$, where $a \neq 0$, $a = 1$ on $\partial\Omega$. Once we prove the stability estimate for $g$ as above, the problem then is reduced to integral geometry problems of recovering the form $\sum b_j dx^j$ (up to $d\phi$ with $\phi|_{\partial\Omega} = 0$) from integrals along geodesics and to that of recovering $q$ from weighted integrals along geodesics, where the weight depends on $b$. Stability estimates for the first one are provided in [BG], see also [StU1]. The second one is injective with natural stability estimates for generic $(g, b)$, as follows from the analysis in [StU2]. Uniqueness of the recovery of $P$, up to the obstructions above, was proven in [KL] without restrictions on $g, b, q$.

The plan of the paper is the following. In section 3 we first prove a Hölder stability estimate at the boundary. We show that if the DN maps of two metrics are close, then their derivatives on the boundary are close, in boundary normal coordinates. This follows essentially from the fact that, away from the glancing manifold, $\Lambda_g$ is locally a pseudo-differential operator, and the normal derivatives of $g$ can be recursively reconstructed in an explicit way from its full symbol, see [SyU]. In section 4, we prove interior stability, i.e., we prove the main result. To this end we prove first a Hölder type of stability estimate that proves that if the DN maps of two metrics are close, then their boundary distance functions are close, too, and then we apply the results in [StU2].

## 2 Preliminaries

**Definition 1** *We say that the Riemannian metric $g$ is simple in $\Omega$, if $\partial\Omega$ is strictly convex w.r.t. $g$, and for any $x \in \bar{\Omega}$, the exponential map $\exp_x : \exp_x^{-1}(\bar{\Omega}) \to \bar{\Omega}$ is a diffeomorphism.*

Note that a small $C^k(\bar{\Omega})$, $k \gg 1$, perturbation of a simple metric in $\Omega$ is also simple. Next, if $g$ is simple, one can extend $g$ in a strictly convex neighborhood $\Omega_1 \supset \bar{\Omega}$ as a simple metric in $\Omega_1$.

The geodesic $X$-ray transform $I_g$ of 2-tensors, which is a linearization of the boundary rigidity problem, is defined as
$$I_g f(\gamma) = \int f_{ij}(\gamma(t)) \dot{\gamma}^i(t) \dot{\gamma}^j(t) \, dt,$$
where $f_{ij}$ is a symmetric tensor. It is known that $I_g d^s v = 0$ for any vector field $v$ with $v = 0$ on $\partial\Omega$. Here $d^s$ is the symmetric differential defined by $[d^s v]_{ij} = \frac{1}{2}(\nabla_i v_j + \nabla_j v_i)$, and $\nabla_i$ are the covariant derivatives.

**Definition 2** *We say that $I_g$ is s-injective in $\Omega$, if $I_g f = 0$ and $f \in L^2(\Omega)$ imply $f = d^s v$ with some vector field $v \in H_0^1(\Omega)$.*

**Definition 3** *Given $k \geq 2$, define $\mathcal{G}^k = \mathcal{G}^k(\Omega)$ as the set of all simple $C^k(\bar{\Omega})$ metrics in $\Omega$ for which the map $I_g$ is s-injective.*

By [StU2], for $k \gg 1$, $\mathcal{G}^k$ is open and dense subset of all simple $C^k(\bar{\Omega})$ metrics, and in particular, all real analytic simple metrics belong to $\mathcal{G}^k$. All metrics with small enough bound on the curvature, and in particular all negatively curved metrics belong to $\mathcal{G}^k$, see [Sh] and the references there. We have local uniqueness for the boundary rigidity problem near metrics in $\mathcal{G}^k$, global uniqueness for pair of metrics in an open and dense set $U \subset \mathcal{G}^k \times \mathcal{G}^k$, and a conditional stability estimate of Hölder type, see [StU2] and (34).

To simplify the notation, we denote

$$\|\cdot\|_* = \|\cdot\|_{H_0^1([0,\varepsilon]\times\partial\Omega)\to L^2([0,T]\times\partial\Omega)}.$$

More precisely, if $\Lambda$ is the DN map defined above, and $f \in H_0^1([0,\varepsilon]\times\partial\Omega)$, then $\|\Lambda\|_*$ is defined as the supremum of $\|\Lambda f\|_{H^1([0,T]\times\partial\Omega)}$ over all $f$ as above with $\|f\|_{H^1([0,\varepsilon]\times\partial\Omega)} = 1$. The correctness of this definition is justified by the following. One can extend any such $f$ as zero for $t > \varepsilon$ and the so extended $f$ will be in $H^1([0,T]\times\partial\Omega)$ with $f|_{t=0} = 0$. We can use standard estimates for mixed hyperbolic problems, see [CP], to show that $\|\Lambda\|_* < \infty$, because it can be estimated by the same norm with $\varepsilon = T$.

## 3 Stability at the boundary

We will prove first stability at the boundary. The arguments here are close to those in [StU3, Prop. 5.1] and [SyU].

Fix a simple metric $g_0 \in C^k$, $k \gg 1$. Extend $g_0$ as a simple metric in some $\Omega_1 \supset \bar{\Omega}$. Let $g$, $\tilde{g}$ be two metrics that will play the role of $g_1$, $g_2$ in Theorem 1 with some $A > 0$ and $\varepsilon_0 \ll 1$, i.e., we have

$$\|g\|_{C^k(\bar{\Omega})} + \|\tilde{g}\|_{C^k(\bar{\Omega})} \leq M, \quad \|g - g_0\|_{C(\bar{\Omega})} + \|\tilde{g} - g_0\|_{C(\bar{\Omega})} \leq \varepsilon_0. \tag{3}$$

The first condition above is a typical compactness condition. Using the interpolation estimate [Tri]

$$\|f\|_{C^t(\bar{\Omega})} \leq C\|f\|_{C^{t_1}(\bar{\Omega})}^{1-\theta}\|f\|_{C^{t_2}(\bar{\Omega})}^{\theta}, \quad t = (1-\theta)t_1 + \theta t_2, \tag{4}$$

where $0 < \theta < 1$, $t_1 \geq 0$, $t_2 \geq 0$, one gets that $\|g - g_0\|_{C^t(\bar{\Omega})} \leq C(M)\varepsilon_0^{(k-t)/k}$ for each $t \geq 0$, if $k > t$; the same is true for $\tilde{g}$. For our purposes, it is enough to apply (4) with $t$, $t_1$ and $t_2$ integers only, then (4) easily extends to compact manifolds with or without boundary. Set

$$\delta = \|\Lambda - \tilde{\Lambda}\|_*. \tag{5}$$

Here and below, a tilde above an object indicates that it is associated with $\tilde{g}$. If there is no tilde, it is related to $g$.

We need here a highly oscillating solution asymptotically supported near a single geodesic transversal to $\partial\Omega$. We need to work only locally near a fixed point $x_0 \in \partial\Omega$, and let $(x', x^n)$ be boundary normal coordinates near $x_0$. Let $\lambda > 0$ be a large parameter. Fix $t_0$ such that $0 < t_0 < \varepsilon/10$, and let $\chi \in C_0^\infty(\mathbf{R}_+ \times \partial\Omega)$ be supported in a small enough neighborhood of $(t_0, x_0)$ of radius not exceeding $\varepsilon/100$ and equal to 1 in a smaller neighborhood of this point. We define $u$ as the solution to (1) with

$$f = e^{i\lambda(t-\phi(x,\omega))}\chi(t,x). \tag{6}$$

One can get an asymptotic expansion of $u$ near $(t_0, x_0)$ by looking for $u$ of the form

$$u = e^{i\lambda(t-\phi(x,\omega))}\sum_{j=0}^{N}\lambda^{-j}A_j(t,x,\omega) + O(\lambda^{-N-1}), \tag{7}$$

where $N \gg 0$ is fixed, $g^{ij}(x_0)\omega_i\omega_j = 1$, $g^{ij}(x_0)\nu_j(x_0)\omega_j < 0$. The phase function solves the eikonal equation

$$\sum_{i,j=1}^{n} g^{ij} \frac{\partial \phi}{\partial x^i} \frac{\partial \phi}{\partial x^j} = 1, \quad \phi|_{\partial \Omega} = x \cdot \omega \tag{8}$$

with the extra condition $\partial \phi/\partial \nu|_{\partial \Omega} < 0$. It is uniquely solvable near $x_0$. In our coordinates, the metric $g$ satisfies $g_{in} = \delta_{in}$, $g^{in} = \delta^{in}$ for $i = 1, \ldots, n$, and $\partial/\partial \nu = -\partial/\partial x^n$. By (8), $\partial \phi/\partial x^n|_{x^n=0} = \omega_n > 0$. The principal part $A_0$ of the amplitude solves near $(t_0, x_0)$ the transport equation (see [SyU])

$$LA_0 = 0, \quad A_0|_{x^n=0} = \chi, \tag{9}$$

and the lower order terms solve

$$iLA_j = (\partial_t^2 - \Delta_g)A_{j-1}, \quad A_j|_{x^n=0} = 0, \quad j \geq 1, \tag{10}$$

where

$$L = 2\partial_t + 2\frac{\partial \phi}{\partial x^n}\frac{\partial}{\partial x^n} + 2\sum_{\alpha,\beta=1}^{n-1} g^{\alpha\beta}\frac{\partial \phi}{\partial x^\alpha}\frac{\partial}{\partial x^\beta} + \Delta_g \phi.$$

The construction of $u$ (see also next section) guarantees that $A_j$, $j = 1, \ldots, N$ are supported in a small neighborhood, depending on the size of supp $\chi$, of the characteristic issued from $(t_0, x_0)$ in the (co-)direction $(1, \omega)$. Therefore, the term $\tilde{u} := e^{i\lambda(t-\phi)} \sum \lambda^{-j} A_j$ in (7) satisfies the zero initial conditions in (1). Moreover, $\tilde{u}$ satisfies the boundary condition $\tilde{u} = f$ with $f$ as in (6), provided that $T$ in (1) is such that $0 < T - t_0$ is small enough. Write $u = \tilde{u} + w$. Then $w = w_t = 0$ for $t = 0$, and $w|_{(0,T)\times \partial \Omega} = 0$ with $T$ as above, and $(\partial_t^2 - \Delta_g)w = O(\lambda^{-N})$. Using standard hyperbolic estimates and Sobolev embedding estimates, one can show that $w = O(\lambda^{-(N-k)})$ in $C^1$, where $k > 0$ depends on $n$ only. We then replace $N$ by $N + k$, and this proves (7) with the estimate of the remainder in the $C^1$ norm. We emphasize that it is important that $T - t_0$ is small enough so that the wave does not meet $\partial \Omega$ again (if it does, we need to reflect it off the boundary, as in next section).

Let $\Psi$, $\tilde{\Psi}$ be two local diffeomorphisms mapping the original coordinates near $x_0$ into boundary normal coordinates $(x', x^n)$ near $(0, 0)$, corresponding to $g$, $\tilde{g}$, respectively. Let $h = \Psi_* g$, $\tilde{h} = \tilde{\Psi}_* \tilde{g}$, and $\varphi = \Psi^{-1}\tilde{\Psi}$. Using a partition of unity, one can extend $\varphi$ in a small neighborhood of $\partial \Omega$.

**Theorem 2** *For any $\mu < 1$, $m \geq 0$, there exists $k > 0$, such that for any $A > 0$, if $\|g_j\|_{C^k(\bar{\Omega})} \leq A$, $j = 1, 2$, then $\exists C > 0$, such that for some diffeomorphism $\varphi$ fixing the boundary,*

$$\sup_{x \in \partial \Omega, |\gamma| \leq m} |\partial^\gamma (g - \varphi_* \tilde{g})| \leq C\delta^{\mu/2^m}.$$

*Proof.* We follow closely [SyU], where it is proven the $\Lambda$ recovers the Taylor series of $g$ at $\partial \Omega$ (i.e., that is the $\delta = 0$ case). Denote $\varphi_* \tilde{g}$ by $\tilde{g}$ again, and work in normal boundary coordinates, the same for both metrics. Observe first that in those coordinates, in a neighborhood $(t_0 - \varepsilon_1, t_0 + \varepsilon_1) \times V \subset \mathbf{R}_+ \times \partial \Omega$ of $t = t_0$, $(x', x^n) = (0, 0)$, where $\chi = 1$, we have

$$\Lambda(u|_{(0,T)\times \partial \Omega}) = e^{i\lambda(t-x\cdot\omega)}\left(i\lambda \frac{\partial \phi}{\partial x^n} - \sum_{j=0}^{N} \lambda^{-j}\frac{\partial}{\partial x^n} A_j\right) + O(\lambda^{-N-1}), \tag{11}$$

and similarly for $\tilde{\Lambda}\tilde{u}$. Therefore,

$$\left\|\frac{\partial \phi}{\partial x^n} - \frac{\partial \tilde{\phi}}{\partial x^n}\right\|_{L^2(V)} \leq \frac{C}{\lambda}\Big(\delta\|u\|_{H^1([0,\varepsilon]\times\partial\Omega)} \tag{12}$$
$$+ \|\tilde{\Lambda}\|_* \|u - \tilde{u}\|_{H^1([0,\varepsilon]\times\partial\Omega)}\Big) + \frac{C}{\lambda}.$$

Notice that $\|u - \tilde{u}\|_{H^1([0,\varepsilon]\times\partial\Omega)} \leq C\lambda^{-N}$, as $\lambda \to \infty$, where $C$ is uniform, if $g$ belongs to a fixed ball in $C^k$ with $k \gg 1$. On the other hand, $\|u\|_{H^1([0,\varepsilon]\times\partial\Omega)} \leq C\lambda$ with a similar $C$. Take the limit $\lambda \to \infty$ above to get

$$\left\|\frac{\partial \phi}{\partial x^n} - \frac{\partial \tilde{\phi}}{\partial x^n}\right\|_{L^2(V)} \leq C\delta.$$

By the eikonal equation (8), in $V \subset \partial\Omega$, we have

$$\frac{\partial \phi}{\partial x^n} = \Big(1 - \sum_{\alpha,\beta=1}^{n-1} g^{\alpha\beta}\omega_\alpha\omega_\beta\Big)^{\frac{1}{2}}, \tag{13}$$

and similarly for $\partial\tilde{\phi}/\partial x^n$. Choosing various $\omega$'s, not tangent to $\partial\Omega$, we prove that $\|g - \tilde{g}\|_{L^2(V)} \leq C\delta$. By a partition of unity argument, this is true on the whole $\partial\Omega$. Using interpolation estimates in Sobolev spaces and Sobolev embedding theorems, we get for any $m \geq 0$ and $\mu < 1$ that

$$\|g - \tilde{g}\|_{C^m(\partial\Omega)} \leq C\delta^\mu, \tag{14}$$

provided that $k \gg 1$, see (3).

To estimate the difference of the first normal derivatives of $g$ and $\tilde{g}$, we use (11) again. As in (12), we have

$$\left\|\frac{\partial A_0}{\partial x^n} - \frac{\partial \tilde{A}_0}{\partial x^n}\right\|_{L^2(V)} \leq C\left(\lambda\delta + \delta + \lambda^{-1}\right). \tag{15}$$

The r.h.s. above is minimized when $\lambda = \delta^{-1/2}$, thus

$$\left\|\frac{\partial A_0}{\partial x^n} - \frac{\partial \tilde{A}_0}{\partial x^n}\right\|_{L^2(V)} \leq C\delta^{1/2}.$$

The transport equation (9) implies that on $x^n = 0$, and on $\chi = 1$, we have

$$\frac{\partial A_0}{\partial x^n} = -\frac{1}{2\omega_n}\frac{1}{\sqrt{\det g}}\frac{\partial}{\partial x^n}\sqrt{\det g}\frac{\partial}{\partial x^n}\phi + R,$$

where $R$ involves tangential derivatives of $g$ only. Therefore,

$$\frac{1}{\sqrt{\det g}}\frac{\partial}{\partial x^n}\sqrt{\det g}\frac{\partial}{\partial x^n}\phi - \frac{1}{\sqrt{\det \tilde{g}}}\frac{\partial}{\partial x^n}\sqrt{\det \tilde{g}}\frac{\partial}{\partial x^n}\tilde{\phi} = O(\delta^{1/2})$$

in $L^2(V)$. By (13), (14),

$$\frac{\partial}{\partial x^n}\sqrt{\det g}\Big(1 - \sum_{\alpha,\beta=1}^{n-1} g^{\alpha\beta}\omega_\alpha\omega_\beta\Big)^{\frac{1}{2}} - \frac{\partial}{\partial x^n}\sqrt{\det \tilde{g}}\Big(1 - \sum_{\alpha,\beta=1}^{n-1} \tilde{g}^{\alpha\beta}\omega_\alpha\omega_\beta\Big)^{\frac{1}{2}} = O(\delta^{1/2})$$

for all $\omega$'s as above. Set $\omega' = 0$ first to estimate the normal derivative of $\det g - \det \tilde{g}$. Choosing finite number of $\omega$'s, we estimate $g^{\alpha\beta} - \tilde{g}^{\alpha\beta}$ for each $\alpha, \beta$ as well. Therefore,

$$\left\|\frac{\partial}{\partial x^n}(g - \tilde{g})\right\|_{C^m(\partial\Omega)} \leq C\delta^{\mu/2} \tag{16}$$

for any $m \geq 0$ and $\mu < 1$ as long as $k \gg 1$.

To estimate the difference of the second normal derivatives of $g$ and $\tilde{g}$, we argue as above. First, we show that, similarly to (15),

$$\left\|\frac{\partial A_1}{\partial x^n} - \frac{\partial \tilde{A}_1}{\partial x^n}\right\|_{L^2(V)} \leq C\left(\lambda^2\delta + \lambda\delta^{1/2} + \lambda^{-1}\right). \tag{17}$$

Choose $\lambda = \delta^{-1/4}$ to get that the r.h.s. above is $O(\delta^{1/4})$. The transport equation (10), including the initial condition $A_1 = 0$ for $x^n = 0$ imply that for $x^n = 0$ and on $\chi = 1$, we have

$$\frac{\partial A_1}{\partial x^n} = -\frac{1}{2\omega_n}\frac{1}{\sqrt{\det g}}\frac{\partial}{\partial x^n}\sqrt{\det g}\frac{\partial}{\partial x^n}A_0.$$

This implies $\partial_{x^n}^2(A_0 - \tilde{A}_0) = O(\delta^{1/4})$ in $L^2(V)$, and therefore the same estimate holds for $\partial_{x^n}^2(\phi - \tilde{\phi})$. This allows us to estimate $\partial_{x^n}^2(g - \tilde{g})$ in the same way. Proceeding by induction, we prove the theorem. $\square$

## 4 Interior Stability

In this section, we prove Theorem 1. The proof is based on the following.

**Proposition 1** *Fix $M > 0$, $\varepsilon_0 > 0$ and let $g$, $\tilde{g}$ be two simple metrics satisfying (3). Then*

$$\|\rho - \tilde{\rho}\|_{C(\partial\Omega \times \partial\Omega)} \leq C\|\Lambda - \tilde{\Lambda}\|_*^{\mu}, \quad \forall x, y \in \partial\Omega$$

*with some $0 < \mu < 1$ depending on $n$ only.*

*Proof.* Recall (5) that $\delta = \|\Lambda - \tilde{\Lambda}\|_*$. It is enough to prove the proposition for $\delta \ll 1$. By Theorem 2, one can assume that for any $m \geq 0$, there exists $\mu > 0$, $k > 0$, such that

$$\sup_{x \in \partial\Omega, |\gamma| \leq m} |\partial^\gamma(g - \tilde{g})| \leq C\delta^\mu. \tag{18}$$

It is convenient to replace the two metrics $g$ and $\tilde{g}$ by two new ones, $g_1$ and $\tilde{g}_1$, equal in a $\delta$-dependent neighborhood of $\partial\Omega$. Let $\chi \in C^\infty(\mathbf{R})$, $\chi(t) = 1$ for $t < 1$, and $\chi(t) = 0$ for $t > 2$. Let $M > 0$ be a large parameter that will be specified below. Set

$$\tilde{g}_1 = \tilde{g} + \chi\left(\delta^{-1/M}\rho(x, \partial\Omega)\right)(g - \tilde{g}), \quad g_1 = g.$$

Using the finite Taylor expansion of $g - \tilde{g}$ of large enough order and (18), we see that (see also [StU2, sec. 7])

$$\|\tilde{g}_1 - \tilde{g}\|_{C^m(\bar{\Omega})} \leq C\delta^{\mu - m/M}.$$

Choose $M = 2m/\mu$. In particular, the estimate above implies that $\tilde{g}_1$ is also simple for $\delta \ll 1$. Without loss of generality we can assume that (3) is still true for $g_1$ and $\tilde{g}_1$. Moreover, one has

$$|\rho_{\tilde{g}}(x, y) - \rho_{\tilde{g}_1}(x, y)| \leq C\delta^{\mu/2}, \quad \forall x, y. \tag{19}$$

We extend $g_1$, $\tilde{g}_1$ in the same way as simple metrics in a neighborhood $\Omega_1 \supset\supset \Omega$. The advantage we have now is that

$$g_1 = \tilde{g}_1 \quad \text{for } -1/C \leq x^n \leq \delta^{1/M}, \tag{20}$$

where $x^n$ is a boundary normal coordinate as above (the same for both metrics).

Next, we will construct an oscillating solution related to $g_1$, similar to the one used in Section 3. To simplify the notation, the objects below related to $g_1$ are without tildes (and without the subscript 1) and those related to $\tilde{g}_1$ have tildes above them (and again, without the subscript 1). Fix $x_0$ and $y_0$ on $\partial\Omega$. Assume that

$$\rho(x_0, y_0) \geq \delta^{\mu/2}. \tag{21}$$

We want to show that $\rho^2(x_0, y_0) - \tilde{\rho}^2(x_0, y_0) = O(\delta^{\mu'})$ with some $0 < \mu' < 1$ under the condition (21) and then to show that this is uniform w.r.t. $x_0$, $y_0$ as in (21). Then we apply the following argument: if a smooth function $f$ on a compact Riemannian manifold with uniformly bounded $C^1$ norm satisfies $f(x) = O(\delta^{\mu'})$ outside a set $W$ with diameter $O(\delta^{\mu''})$, then $f = O(\delta^{\mu'}) + O(\delta^{\mu''})$ as can be easily seen by integrating the derivative of $f$ along curves connecting an arbitrary point in $W$ with a point outside $W$.

All constants below will be uniform w.r.t. $g$ and $\tilde{g}$ satisfying (3) and in particular, independent of the choice of $x_0$ and $y_0$.

Consider the geodesic connecting $x_0$ and $y_0$, extended from $\Omega$ to $\Omega_1$. Let $z_0 \in \Omega_1 \setminus \Omega$ be a point on this geodesic such that the geodesic segment $[z_0, x_0]$ is in $\Omega_1 \setminus \Omega$. We assume that $\rho(z_0, x_0) > 1/C$ with $C > 0$ fixed. Set $\phi(x) = \rho(x, z_0)$. Then, by the simplicity assumption, since $z_0 \in \Omega_1$, we have that $\phi \in C^{k-1}(\bar{\Omega})$, and $\phi$ solves the eikonal equation

$$\sum_{i,j=1}^{n} g^{ij} \frac{\partial \phi}{\partial x^i} \frac{\partial \phi}{\partial x^j} = 1. \tag{22}$$

Then we construct a solution $u$ of (1) of the form

$$u = e^{i\lambda(\phi(x)-t)} \left( A(t, x) + v(t, x; \lambda) \right), \tag{23}$$

where

$$\|v(t, \cdot; \lambda)\|_{C^2} \leq \frac{C}{\lambda}. \tag{24}$$

The construction of $u$ is the same as that in the preceeding section, except that the phase function has different initial condition and we want to solve it all the way to the opposite side of $\partial\Omega$. The principal part $A$ of the amplitude solves the transport equation (9). This is an ODE along the geodesics issued from points in $\text{supp } A \cap (\mathbf{R} \times \partial\Omega)$ in (co-)directions $\nabla\phi$. Let

$$U = \left\{ (t, x) \in \mathbf{R}_+ \times \partial\Omega; \ |t - t_0| + \rho(x, x_0) < \delta^{\mu/2}/C \right\}, \tag{25}$$

where $0 < t_0 \ll 1$ is fixed, and $C \gg 1$ will be specified later. Choose a cut-off function $0 \leq \chi \leq 1$, $\chi \in C_0^\infty(\mathbf{R}_+ \times \partial\Omega)$ such that $\text{supp } \chi \subset U$, and $\chi = 1$ in a set defined as $U$ but with $C$ replaced by $2C$. One can arrange that $|\partial_t \chi| + |\nabla \chi| \leq C\delta^{-\mu/2}$. Then we solve the transport equation (9) with initial condition

$$A|_{\mathbf{R}_+ \times \partial\Omega} = \chi.$$

The solution is supported in a neighborhood of the geodesic connecting $x_0$ and $y_0$ of size $O(\delta^{\mu/2})$, and can be extended all the way to some neighborhood of $y_0$ by the simplicity assumption. If $C$ in (25) is large enough, then $\text{supp } A \cap (\mathbf{R}_+ \times \partial\Omega)$ consists of two disjoint components near $x_0$ and $y_0$ respectively, one of them being $\text{supp } \chi \subset U$. The other one, let us call it $V$, is the image of $U$ under translation by all geodesics issued from $z_0$, passing through $U$. Because of the strict convexity of $\partial\Omega$, each component is of size $O(\delta^{\mu/2})$, at a distance bounded by below by the same quantity. Denote by $B(y, r)$ the ball centered at $y$ with radius $r$. Then $V$ contains the set $V_0 = V \cap B(y_0, \delta^{\mu/2}/C_0)$, such that on $V_0$, we have $A \geq 1/C$. Above, $C_0$ is chosen so that $V_0$ is contained in the translation of the set $\{\chi = 1\} \subset U$ under the geodesics issued from $z_0$.

Then one gets a solution with the required properties except that $u$ does not necessarily vanish in $\Omega$ for $t < 0$ small enough but the principal part $e^{i\lambda(\phi-t)}A$ does. To justify (24) in the $C^2$ norm, we construct lower order terms, similarly to (11), up to order $N$ large enough so that after applying standard a priori estimates and Sobolev embedding estimates, we get the estimate in the $C^2$ norm.

We reflect $u$ off $\partial\Omega$ at $V$ by setting $h = u|_V$ and solving (1) with boundary data $-h$. Let us call this solution $v$. Then $v$ has the form (23) as well, with a different amplitude $B$ instead of $A$ and a phase function $\varphi$. The phase function $\varphi$ still solves the eikonal equation (22) with boundary condition $\varphi|_{\partial\Omega} = \phi$ and is the unique solution with gradient pointing towards the interior of $\Omega$ (the other solution is $\phi$). For $B$ we have $B|_V = -A|_V$. Then $w := u + v$ vanishes on $V$ modulo $O(\lambda^{-1})$, and is well defined for $0 < t < \tau := \rho(x_0, y_0) + \delta^{\mu/2}/C$, $C \gg 1$. We then extend $w$ by imposing zero boundary conditions for $\tau \leq t \leq T$. Clearly, the requirements on $T$ and $\varepsilon$ imply that $T > \tau$, if $\delta \ll 1$. As above, to justify the estimate on the remainder, we need to construct the lower order terms up to some order, as well.

We claim that
$$\frac{1}{\lambda}\|w\|_{H^1([0,T]\times\partial\Omega)} \leq C + \frac{C(\delta)}{\lambda}, \tag{26}$$

where the first constant $C$ is independent of $\delta$ (but it depends on $M$ and $\varepsilon_0$ in (3) as mentioned above). To prove this, we only need to estimate $A\nabla\phi$ on $U$, and $A\nabla\phi$ and $B\nabla\varphi$ on $V$. By the definition of $\phi$, we have $|\nabla\phi| \leq C$ on $U$. Next, we have $|A\nabla\phi| = |B\nabla\varphi| \leq C$ on $V$ as well.

Next, we construct a similar solution $\tilde{w}$ related to $\tilde{g}_1$. We construct first a phase function $\tilde{\phi}$ as $\tilde{\phi}(x) = \tilde{\rho}(x, z_0)$. It solves the eikonal equation
$$\sum_{i,j=1}^n \tilde{g}_1^{ij}\frac{\partial\tilde{\phi}}{\partial x^i}\frac{\partial\tilde{\phi}}{\partial x^j} = 1, \quad \tilde{\phi}|_U = \phi. \tag{27}$$

The latter equality follows from (20). The other properties of $\tilde{w}$ are similar to those of $w$. Let $\tilde{V}$ be defined as above, but associated to $\tilde{g}_1$.

On $V$, we have
$$\Lambda(w|_{(0,T)\times\partial\Omega}) = -2i\lambda e^{i\lambda(\phi-t)}\frac{\partial\phi}{\partial x^n}A + O_\delta(1),$$

and similarly for $\tilde{\Lambda}\tilde{w}$. Notice that
$$\left\|\frac{\partial\phi}{\partial x^n}A\right\|_{L^2(V)} \geq \delta^{n\mu/4}/C, \tag{28}$$

because $\text{area}(V_0) \geq \delta^{n\mu/2}/C$.

If $V \cap \tilde{V} = \emptyset$, then we get by (28),
$$\|\Lambda(w|_{(0,T)\times\partial\Omega}) - \tilde{\Lambda}(\tilde{w}|_{(0,T)\times\partial\Omega})\| \geq \lambda\delta^{n\mu/4}/C - C(\delta). \tag{29}$$

On the other hand,

$$\|\Lambda(w|_{(0,T)\times\partial\Omega}) - \tilde{\Lambda}(\tilde{w}|_{(0,T)\times\partial\Omega})\| \leq \delta\|w\|_{H^1([0,T]\times\partial\Omega)} + C\|w - \tilde{w}\|_{H^1([0,T]\times\partial\Omega)}. \tag{30}$$

Notice that $\|w\|_{H^1([0,T]\times\partial\Omega)}$ is $O_\delta(1)$ outside $U$, as $\lambda \to \infty$. Restricted to $U$, we get $O(\lambda\delta^{-\mu/2+n\mu/4})$. Thus,

$$\|w\|_{H^1([0,T]\times\partial\Omega)} \leq C\lambda\delta^{-\mu/2+n\mu/4} + C(\delta). \tag{31}$$

On the other hand,

$$\|w - \tilde{w}\|_{H^1([0,T]\times\partial\Omega)} \leq C(\delta). \tag{32}$$

Combine the inequalities (29), (30), (31) and (32), to get

$$\lambda\delta^{n\mu/4} \leq C\lambda\delta^{1-\mu/2+n\mu/4} + C(\delta).$$

Divide by $\lambda$ and take the limit $\lambda \to \infty$ to get a contradiction.

The contradiction above shows that $V$ and $\tilde{V}$ do intersect provided that $\delta \ll 1$. Therefore, there exists $q \in V \cap \tilde{V} \subset \partial\Omega$, and $p, \tilde{p} \in U$, such that $\rho(p,q) = \tilde{\rho}(\tilde{p}, q)$. Since the diameters of $U$, $V$, and $\tilde{V}$ are $O(\delta^{\mu/2})$, we get that

$$|\rho(x_0, y_0) - \tilde{\rho}(x_0, y_0)| \leq C\delta^{\mu/2}. \tag{33}$$

Recall now that by our notation convention, $\rho = \rho_{g_1}$, $\tilde{\rho} = \rho_{\tilde{g}_1}$ above. Combine this with (19) and the argument following (21), to complete the proof of the proposition. □

*Proof of Theorem 1:* The proof follows directly by combining Proposition 1 and Theorem 4 in [StU2]. Indeed, under the assumptions of Theorem 1, it was shown in [StU2, Theorem 5], that for any $\mu < 1$,

$$\|g_2 - \psi^* g_1\|_{C^2(\bar{\Omega})} \leq C(A)\|\rho_{g_1} - \rho_{g_2}\|^\mu_{C(\partial\Omega\times\partial\Omega)} \tag{34}$$

as long as $k \gg 1$. Apply Proposition 1 to complete the proof. □